\documentclass[12pt]{article}

\usepackage{graphics}
\usepackage{latexsym}
\usepackage{amsmath}
\usepackage{amsfonts}
\usepackage{booktabs}
\usepackage{stackrel}
\usepackage{color}
\usepackage{appendix}

\def\ds{\displaystyle}
\def\dt{\frac{\partial }{\partial t}}
\def\dx{\frac{\partial }{\partial x}}

\def\hat{\widehat}

\def\bar{\overline}

\def\bar{\overline}
\def\mbf{\mathbf}

\def\hat{\widehat}

\def\bar{\overline}
\def\Re{{\rm Re}}
\def\Im{{\rm Im}}
\def\Re{{\rm Re}}
\def\Im{{\rm Im}}
\newcommand{\R}{{\mathbb R}}
\newcommand{\C}{{\mathbb C}}

\newtheorem{theorem}{Theorem}[section]

\def\phih{\hat{\phi}}

\usepackage{amssymb}
\usepackage{amsmath}
\usepackage{graphicx}

\begin{document}

\title{Transformed implicit-explicit DIMSIMs with strong stability preserving explicit part}

\author{G. Izzo%
\thanks{Dipartimento di Matematica e Applicazioni,
Universit\`a di Napoli Federico II, 
80126 Napoli, Italy,
\mbox{e-mail}: giuseppe.izzo@unina.it.
Member of the INdAM Research group GNCS.
} \
and Z. Jackiewicz,%
\thanks{Department of Mathematics, Arizona State University,
Tempe, Arizona 85287, and AGH University of Science and Technology, Krak\'ow, Poland,
\mbox{e-mail}: jackiewicz@asu.edu.} }

\date{\today}

%
%
%
%

\maketitle

\textbf{Abstract.} 
For many systems of differential equations modeling problems 
in science and 
engineering, there are often natural splittings of the right hand side into two
parts, one of which is non-stiff or mildly stiff, and 
the other part is stiff.
Such systems can be efficiently treated by a class of 
implicit-explicit (IMEX) 
diagonally implicit multistage integration methods (DIMSIMs),
where the stiff part is integrated by an implicit formula, and the 
non-stiff part is
integrated by an explicit formula. 
We will construct methods where the explicit part has strong stability preserving
(SSP) property, and the implicit part of the method is $A$-, or $L$-stable.
We will also investigate stability of these methods 
when the implicit and explicit parts interact with each other. 
To be more precise,
we will monitor the size of the region of absolute stability
of the IMEX scheme, assuming that the 
implicit part of the method is $A$-, or $L$-stable. Finally we furnish 
examples of SSP IMEX DIMSIMs up to the order four with good stability properties.
\vspace{1mm}


\vspace{2mm}

\newpage

\setcounter{equation}{0}
\setcounter{figure}{0}
\setcounter{table}{0}

\section{Introduction} \label{sec1}
Many practical problems in science and engineering are modeled by large systems
of ordinary differential equations (ODEs) which arise from discretization in space
of partial differential equations (PDEs) by finite difference methods, 
finite elements or finite volume methods, or pseudospectral methods.
For such systems there are often natural splittings of the right hand sides 
of the differential systems into two parts, one of which is non-stiff or
mildly stiff, and suitable for explicit time integration, and the other part is stiff, 
and suitable for implicit time integration. Such systems can be written in the form
\begin{equation} \label{eq1.1}
\begin{array}{ll}
y'(t)=f\big(y(t)\big)+g\big(y(t)\big), & t\in[t_0,T], \\
y(t_0)=y_0\in\R^m,
\end{array}
\end{equation}
$f:\R^m\rightarrow \R^m$,
$g:\R^m\rightarrow \R^m$,
where $f(y)$ represents the non-stiff processes, for example advection, 
and $g(y)$ represents stiff processes, for example diffusion or chemical reaction,
in semidiscretization of advection-diffusion-reaction equations \cite{hv03}.

In this paper we will analyze methods, where the non-stiff part $f(y)$ is
treated by the explicit general linear method (GLM)
and the stiff part $g(y)$ by the implicit GLM, with the same abscissa vector 
$\mbf{c}=[c_1,\ldots,c_s]^T\in\R^s$, and the coefficients
$$
\left[
\begin{array}{c|c}
\mbf{A} & \mbf{U} \\
\hline
\mbf{B} & \mbf{V}
\end{array}
\right]\in\R^{(s+r)\times(s+r)},
\quad
\left[
\begin{array}{c|c}
\mbf{A}^{*} & \mbf{U} \\
\hline
\mbf{B^{*}} & \mbf{V}
\end{array}
\right]\in\R^{(s+r)\times(s+r)},
$$
We assume that both methods have the same coefficients matrices
$\mbf{U}$ and $\mbf{V}$, and that
$\mbf{A}$ is strictly lower triangular, 
and $\mbf{A}^{*}$ is lower triangular with the same element
$\lambda>0$ on the diagonal.
Denote the components of $\mbf{A}$, $\mbf{A}^{*}$,
$\mbf{U}$, $\mbf{B}$, $\mbf{B}^{*}$, and $\mbf{V}$
by $a_{ij}$, $a^{*}_{ij}$, $u_{ij}$, $b_{ij}$, $b^{*}_{ij}$, 
and $v_{ij}$.
Then on the uniform grid $t_n=t_0+nh$, $n=0,1,\ldots,N$,
$Nh=T-t_0$, the IMEX GLMs are defined by
\begin{equation} \label{eq1.2}
\begin{array}{l}
Y_i^{[n+1]}=h\ds\sum_{j=1}^{i-1}a_{ij}f\big(Y_j^{[n+1]}\big)
+h\ds\sum_{j=1}^ia^{*}_{ij}g\big(Y_j^{[n+1]}\big)
+\ds\sum_{j=1}^ru_{ij}y_j^{[n]},
\ i=1,2,\ldots,s, \\
y_i^{[n+1]}=h\ds\sum_{j=1}^s\Big(b_{ij}f\big(Y_j^{[n+1]}\big)
+b^{*}_{ij}g\big(Y_j^{[n+1]}\big)\Big)
+\ds\sum_{j=1}^rv_{ij}y_j^{[n]},
\quad i=1,2,\ldots,r,
\end{array}
\end{equation}
$n=0,1,\ldots,N-1$.
Here,
$Y_i^{[n+1]}$ are approximations of stage order $q$ to $y(t_n+c_ih)$,
i.e., 
\begin{equation} \label{eq1.3}
Y_i^{[n+1]}=y(t_n+c_ih)+O(h^{q+1}),
\quad
i=1,2,\ldots,s,
\end{equation}
and $y_i^{[n]}$ are approximations of order $p$ to the 
linear combinations of the derivatives of the solution $y$ at the point $t_n$,
i.e.,
\begin{equation} \label{eq1.4}
y_i^{[n]}=\ds\sum_{k=0}^pq_{ik}h^ky^{(k)}(t_n)+O(h^{p+1}),
\quad
i=1,2,\ldots,r,
\end{equation}
where $y$ is the solution to (\ref{eq1.1}).
These IMEX methods were introduced in \cite{zsb14} and further investigated in \cite{cjsz15}.

Putting
$$
y^{[n+1]}=\left[
\begin{array}{c}
y_1^{[n+1]} \\
\vdots \\
y_r^{[n+1]}
\end{array}
\right],
\quad
y^{[n]}=\left[
\begin{array}{c}
y_1^{[n]} \\
\vdots \\
y_r^{[n]}
\end{array}
\right],
\quad
Y^{[n+1]}=\left[
\begin{array}{c}
Y_1^{[n+1]} \\
\vdots \\
Y_s^{[n+1]}
\end{array}
\right],
$$
$$
f\big(Y^{[n+1]}\big)=\left[
\begin{array}{c}
f\big(Y_1^{[n+1]}\big) \\
\vdots \\
f\big(Y_s^{[n+1]}\big)
\end{array}
\right],
\quad
g\big(Y^{[n+1]}\big)=\left[
\begin{array}{c}
g\big(Y_1^{[n+1]}\big) \\
\vdots \\
g\big(Y_s^{[n+1]}\big)
\end{array}
\right],
$$
the method (\ref{eq1.2}) can be written in a more compact form
\begin{equation} \label{eq1.5}
\begin{array}{l}
Y^{[n+1]}=h(\mbf{A}\otimes \mbf{I})f\big(Y^{[n+1]}\big)+h(\mbf{A}^{*}\otimes \mbf{I})g\big(Y^{[n+1]}\big)
+(\mbf{U}\otimes \mbf{I})y^{[n]}, \\ [3mm]
y^{[n+1]}=h(\mbf{B}\otimes \mbf{I})f\big(Y^{[n+1]}\big)+h(\mbf{B}^{*}\otimes \mbf{I})g\big(Y^{[n+1]}\big)
+(\mbf{V}\otimes \mbf{I})y^{[n]},
\end{array}
\end{equation}
$n=0,1,\ldots,N-1$, $\mbf{I}\in\R^m$, and the relation (\ref{eq1.4}) takes the form
\begin{equation} \label{eq1.6}
y^{[n]}=\ds\sum_{k=0}^p\mbf{q}_kh^ky^{(k)}(t_n)+O(h^{p+1}),
\end{equation}
with the vectors $\mbf{q}_0,\mbf{q}_1,\ldots,\mbf{q}_s$ given by
$$
\mbf{q}_0=\left[
\begin{array}{c}
q_{1,0} \\
\vdots \\
q_{r,0}
\end{array}
\right],
\quad
\mbf{q}_1=\left[
\begin{array}{c}
q_{1,1} \\
\vdots \\
q_{r,1}
\end{array}
\right],
\quad
\ldots,
\quad
\mbf{q}_p=\left[
\begin{array}{c}
q_{1,p} \\
\vdots \\
q_{r,p}
\end{array}
\right].
$$

In this paper we will investigate the class of IMEX
diagonally implicit multistage integration methods (DIMSIMs).
These are schemes with $p=q=r=s$, where the coefficient matrix
$\mbf{U}=\mbf{I}$, and $\mbf{V}$ is  a rank one matrix of the form
$\mbf{V}=\mbf{e}\mbf{v}^T$, $\mbf{e}=[1,\ldots,1]^T\in\R^s$,
$\mbf{v}=[v_1,\ldots,v_s]^T\in\R^s$, with $\mbf{v}^T\mbf{e}=1$.
It was proved in \cite{but93} (see also \cite{jac09}) that for
given $\mbf{A}$, $\mbf{A}^{*}$, and $\mbf{V}$ the explicit 
method and the implicit method has order $p$ and stage order
$q=p$ if the coefficients matrices $\mbf{B}$ and $\mbf{B}^{*}$
are computed from the formulas
\begin{equation} \label{eq1.7}
\mbf{B}=\mbf{B}_0-\mbf{A}\mbf{B}_1-\mbf{V}\mbf{B}_2+\mbf{V}\mbf{A},
\quad
\mbf{B}^{*}=\mbf{B}_0-\mbf{A}^{*}\mbf{B}_1-\mbf{V}\mbf{B}_2+\mbf{V}\mbf{A}^{*},
\end{equation}
where $\mbf{B}_0$, $\mbf{B}_1$, and $\mbf{B}_2$, are $s\times s$ matrices
defined by
$$
\mbf{B}_0=\left[
\ds\frac{\int_0^{1+c_i}\phi_j(x)dx}{\phi_j(c_j)}
\right],
\quad
\mbf{B}_1=\left[
\ds\frac{\phi_j(1+c_i)}{\phi_j(c_j)}
\right],
\quad
\mbf{B}_2=\left[
\ds\frac{\int_0^{c_i}\phi_j(x)dx}{\phi_j(c_j)}
\right],
$$
$i,j=1,2,\ldots,s$, and $\phi_i(x)$ are defined by
$$
\phi_i(x)=\prod_{j=1,j\neq i}^s(x-c_j),
\quad
i=1,2,\ldots,s.
$$
It was also proved in \cite{zsb14} that if the explicit and implicit methods
have order $p$ and stage order $q=p$, then the same is true for the
resulting IMEX scheme defined by (\ref{eq1.5}).

The methods investigated in this paper are also applicable
to the hyperbolic systems with relaxation considered, for example,
in \cite{jin95,pr05}, and they compare favorably with 
IMEX Runge-Kutta (RK) methods for these problems.
In the stiff limit the IMEX RK schemes, considered for example
in \cite{pr05} converge, but their order drops to $p=1$,
while all IMEX DIMSIMs constructed in this paper achieve
the expected order of convergence, and no order reduction 
occurs. This is confirmed in Section~\ref{sec5} by 
numerical experiments on the shallow water equation.

The organization of the remainder of the paper is as follows.
In Section~\ref{sec2} we will review various stability concepts
of explicit, implicit, and IMEX schemes.
In particular, we will recall the definition of strong stability
preserving (SSP) property of explicit methods, 
absolute stability, and definitions of regions of absolute stability,
of explicit, implicit, and the resulting IMEX 
methods. In Section~\ref{sec3} we define transformed IMEX methods. 
In Section~\ref{sec4} we describe the construction of SSP transformed IMEX
schemes of order $p=1$, $2$, $3$, and $4$.
In Section~\ref{sec5} the results of some numerical experiments are 
presented.

\setcounter{equation}{0}
\setcounter{figure}{0}
\setcounter{table}{0}

\section{Stability analysis of IMEX DIMSIMs} \label{sec2}
\subsection{SSP property of the explicit part} \label{sec2.1}
We recall first the concept of SSP property of explicit methods
following the presentation in \cite{gks11}.
To define this property we assume that the discretization of the problem
(\ref{eq1.1}) with $g\equiv 0$, by the forward Euler method
$$
y_{n+1}=y_n+hf(y_n),
\quad
n=0,1,\ldots,N-1,
$$
satisfies the inequality
\begin{equation} \label{eq2.1}
\|y_{n+1}\|\leq \|y_n\|,
\quad
n=0,1,\ldots,N-1,
\end{equation}
in some norm or semi-norm $\|\cdot\|$,
if the time step $h$ is restricted by the condition
\begin{equation} \label{eq2.2}
h\leq h_{FE}. 
\end{equation}
It is then of interest to construct higher order numerical
methods for (\ref{eq1.1}) with $g\equiv 0$, which preserve the property
(\ref{eq2.1}) under the time step restrictions
\begin{equation} \label{eq2.3}
h\leq \mathcal{C}\cdot h_{FE},
\end{equation}
where $\mathcal{C}\geq 0$ is some constant.
Numerical schemes for (\ref{eq1.1}) with
$g\equiv 0$, which preserve the property
(\ref{eq2.1}) under the condition (\ref{eq2.3}) 
are called SSP methods, and the maximal constant $\mathcal{C}$ in
(\ref{eq2.3}) is called SSP coefficient. To compare
numerical methods with different number of stages $s$ we also define,
following \cite{cs10,gks11,kgm11}, the effective SSP coefficient 
$\mathcal{C}_{eff}$ by the relation $\mathcal{C}_{eff}=\mathcal{C}/s$.

The characterization of SSP coefficient for GLMs was discovered by
Spijker \cite{spi07}.
To describe this characterization for GLMs defined by
the abscissa vector $\mbf{c}$ and coefficient matrices
$\mbf{A}$, $\mbf{U}$, $\mbf{B}$, and $\mbf{V}$, consider the relations
\begin{equation} \label{eq2.4}
\begin{array}{ll}
(\mbf{I}+\gamma\mathbf{A})^{-1}\mbf{U}\geq 0, &
\mbf{I}-(\mbf{I}+\gamma\mbf{A})^{-1}\geq 0, \\
\mbf{V}-\gamma\mbf{B}(\mbf{I}+\gamma\mbf{A})^{-1}\mbf{U}\geq 0, &
\gamma\mbf{B}(\mbf{I}+\gamma\mbf{A})^{-1}\geq 0,
\end{array} 
\end{equation}
where $\gamma\geq 0$ is a constant, and where these inequalities
should be interpreted componentwise. Then it was demonstrated by
Izzo and Jackiewicz
\cite{ij15}, using the results by Spijker \cite{spi07}, that the SSP coefficient is given by
\begin{equation} \label{eq2.5}
\mathcal{C}=\mathcal{C}(\mbf{c},\mbf{A},\mbf{U},\mbf{B},\mbf{V})=
\sup\Big\{\gamma\in\R: \ \gamma \ \textrm{satisfies} \ (\ref{eq2.4})\Big\}.
\end{equation}
It follows from this relation that SSP coefficient $\mathcal{C}$ can be computed
by solving the minimization problem
\begin{equation} \label{eq2.6}
F(\gamma):=-\gamma\longrightarrow \min,
\end{equation}
with a very simple objective function
$F(\gamma)=-\gamma$, subject to the nonlinear constrains
(\ref{eq2.4}). This process will be illustrated in Section~\ref{sec4}.

\subsection{Absolute stability of the implicit part} \label{sec2.2}
DIMSIMs investigated in the literature have so-called RK stability
property, i.e., their stability function $p^{*}(w,z)$ assumes the form
$$
p^{*}(w,z)=w^{s-1}(w-R(z)),
$$
where $R(z)$ is an approximation of order $p$ to the exponential
function $\exp(z)$.
However, this stability requirement is quite restrictive and does not,
in general, permit the construction
of IMEX schemes with $A$- or $L$-stable implicit part and SSP explicit part, 
and, for this reason, we do not enforce RK stability of implicit methods in this paper.
But we will still refer to the resulting implicit formulas as DIMSIMs.

In order to construct methods with $A$- or $L$-stable implicit part,
we will apply the the well known 
Schur criterion (\cite{s16}) in combination with the 
maximum principle.
Let us recall that the polynomial
$$
\phi_k(w)= c_k w^k + c_{k-1}w^{k-1} + \cdots + c_1 w + c_0 
$$
where $c_i$ are complex coefficients, with $c_k \neq 0$ and $c_0\neq 0$, 
is said to be a Schur polynomial 
if all of its roots $w_i$, $i = 1, 2, \ldots , k$ , are inside of the unit circle,
that is $|w_i|<1$, for all $i = 1, 2, \ldots , k$.
Define the following two polynomials 
$$
\phih_k(w)= \bar{c}_0 w^k + \bar{c}_{1}w^{k-1} + \cdots + \bar{c}_{k-1} w + \bar{c}_k, 
$$
and
$$
\phi_{k-1}(w)=\frac{1}{w}\left( \phih(0)\phi(w) - \phi(0)\phih(w) \right),
$$
where $\bar{c}_i$ represents the complex conjugate of the coefficient $c_i$, $i=0,1,\ldots,k$,
and let us remark that the polynomial $\phi_{k-1}(w)$ has degree at most $k-1$.
The Schur recursive criterion is based on the following result.

\begin{theorem} (J. Schur \cite{s16})
 The polynomial $\phi_k(w)$ is a Schur polynomial if and only if
$$
\big|\phih_k(0)\big|>\big|\phi_k(0)\big|,
$$ 
and $\phi_{k-1}(w)$ is a Schur polynomial.
\end{theorem}

To analyze stability properties of implicit methods it is convenient to
multiply stability function of these methods by the polynomial factor $(1-\lambda z)^s$.
The resulting stability polynomial $p^{*}(w,z)$ takes then the form
\begin{equation} \label{eq2.7}
p^{*}(w,z)=(1-\lambda z)^sw^s-p_1(z)w^{s-1}+p_2(z)w^{s-2}+\cdots+(-1)^sp_s(z),
\end{equation}
where $p_1(z),p_2(z),\ldots,p_s(z)$ are polynomials of degree less than or equal to $s$.
To construct implicit formulas whose stability polynomial (\ref{eq2.7}) is a Schur polynomial 
in the left half of the complex plane, 
we will force all the roots $w_j = w_j(z)$, $j=1,2,\ldots,r$, of $p^{*}(w,z)$ 
to have no poles for $Re(z)\leq 0$. 
Since these roots are analytic functions of $z$ for $Re(z) \leq 0$, 
they fall inside the unit circle for $Re(z) \leq 0$ if and only if they 
are inside the unit circle for the values of $z$ on the imaginary axis. 
In other words, by the maximum principle (compare \cite{but03}), it follows that
$|w_j(z)|<1$, $j=1,2,\ldots,r$, for all $z\in\C$ with $Re(z)\leq 0$, 
if and only if $|w_j(i y)|<1$, $j=1,2,\ldots,r$, for all $y\in\R$. 

For methods with number of stages $s=r$, the stability polynomial $p^*(w,z)$
has degree $r$ and the conditions given by the recursive Schur criterion 
are the following
$$
\big|\phih_r(0)\big|-\big|\phi_r(0)\big|>0,
\quad
\big|\phih_{r-1}(0)\big|-\big|\phi_{r-1}(0)\big|>0,
\quad
\ldots,
\quad
\big|\phih_{1}(0)\big|-\big|\phi_{1}(0)\big|>0.
$$
Let us define the quantities 
$$
a_k:=\big|\phih_{r-k}(0)\big|-\big|\phi_{r-k}(0)\big|, 
\quad
k=0,1,\dots,r-1.
$$
Each $a_k$ depends on $z$, and when it is evaluated at $z=i y$, $y \in \R$,
it results to be a polynomial in the unknown $y$ with real coefficients, of the form
$$
a_k(i y)=\sum_{j=0}^{r2^k}m_{kj}y^{2j}, 
\quad 
k=0,1,\ldots,r-1,
$$
with $m_{kj}\in\R$ for all $k=0,1,\ldots,r-1$ and $j=1,2,\ldots,r2^k$.
Thus, a sufficient condition to ensure the $A$-stability of the corresponding 
method is to force 
\begin{equation}\label{eq2.8}
m_{kj} \geq 0,
\quad
k=1,2,\ldots,r,
\quad
j=1,2,\ldots,r2^k,
\end{equation}
where for each $k$ at least one $m_{kj}$, $j=0,1,\ldots,r2^k$, has to be strictly positive.

To construct methods which are $L$-stable we have to enforce the condition that
the polynomials $p_1(z),p_2(z),\ldots,p_s(z)$ appearing in (\ref{eq2.7})
have degrees strictly less than $s$.

\subsection{Absolute stability of the IMEX method} \label{sec2.3}
We will discuss next absolute stability.
To analyze absolute stability properties of IMEX GLMs (\ref{eq1.5}) we will
use the test equation
\begin{equation} \label{eq2.9}
y'(t)=\lambda_0y(t)+\lambda_1y(t),
\quad
t\geq 0,
\end{equation}
where $\lambda_0$ and $\lambda_1$ are complex parameters.
Here, $\lambda_0 y(t)$ corresponds to the non-stiff part and 
$\lambda_1 y(t)$ to the stiff part of the system (\ref{eq1.1}).
Applying (\ref{eq1.5}) to (\ref{eq2.9}) and putting 
$z_0=h\lambda_0$, $z_1=h\lambda_1$, we obtain
$$
\begin{array}{l}
Y^{[n+1]}=(z_0\mbf{A}+z_1\mbf{A}^{*})Y^{[n+1]}+\mbf{U}y^{[n]}, \\
y^{[n+1]}=(z_0\mbf{B}+z_1\mbf{B}^{*})Y^{n+1]}+\mbf{V}y^{[n]},
\end{array}
$$
$n=0,1,\ldots$.
Assuming that the matrix $\mbf{I}-z_0\mbf{A}-z_1\mbf{A}^{*}$ is nonsingular,
this is equivalent to the vector recurrence relation
\begin{equation} \label{eq2.10}
y^{[n+1]}=\mbf{M}(z_0,z_1)y^{[n]}, 
\end{equation}
$n=0,1,\ldots$, with the stability matrix $\mbf{M}(z_0,z_1)$ defined by
\begin{equation} \label{eq2.11}
\mbf{M}(z_0,z_1)=\mbf{V}+(z_0\mbf{B}+z_1\mbf{B}^{*})
(\mbf{I}-z_0\mbf{A}-z_1\mbf{A}^{*})^{-1}\mbf{U}.
\end{equation}
We also define the stability function $p(w,z_0,z_1)$ of the IMEX scheme (\ref{eq1.5})
as the stability polynomial of $\mbf{M}(z_0,z_1)$, i.e., 
\begin{equation} \label{eq2.12}
p(w,z_0,z_1)=\det\big(w\mbf{I}-\mbf{M}(z_0,z_1)\big).
\end{equation}
To investigate stability properties of (\ref{eq1.5}) is is usually more convenient
to work with the polynomial $(1-\lambda z_1)^sp(w,z_0,z_1)$, where 
$\lambda$ is the diagonal element of the coefficient matrix $\mbf{A}^{*}$.
This polynomial will be denoted
by the same symbol $p(w,z_0,z_1)$.

We say that the IMEX GLM (\ref{eq1.5}) is stable for given $z_0,z_1\in\C$ if
all the roots $w_i(z_0,z_1)$, $i=1,2,\ldots,r$, of the stability function
$p(w,z_0,z_1)$ are inside of the unit circle. 
In this paper 
we will be mainly interested in IMEX schemes which 
are $A$-stable with respect to the implicit part $z_1\in\C$. 
To investigate such methods 
we consider, similarly as in \cite{cjsz14a,hv03,zsb14}, the sets
$$
\mathcal{S}_{\alpha}=
\Big\{z_0\in \C: \ \textrm{the IMEX GLM is stable for any} \ z_1\in \mathcal{A}_{\alpha}\Big\},
$$
where the set $\mathcal{A}_{\alpha}\subset\C$ is defined by
$$
\mathcal{A}_{\alpha}=\Big\{z\in\C: \ \Re(z)<0 \quad \textrm{and} \quad
|\Im(z)|\leq \tan(\alpha)|\Re(z)|\Big\}.
$$
It follows from the maximum principle that $\mathcal{S}_{\alpha}$ has
a simple representation given by
\begin{equation} \label{eq2.13}
\mathcal{S}_{\alpha}=
\left\{
\begin{array}{lcl}
z_0\in \C\!\!\!\!&:&  \textrm{the IMEX GLM is stable for any} \\
& & z_1=-|y|/\tan(\alpha)+iy, \ y\in\R
\end{array}
\right\}.
\end{equation}
For fixed values of $y\in \R$ we define also the sets
\begin{equation} \label{eq2.14}
\mathcal{S}_{\alpha,y}=
\left\{
\begin{array}{lcl}
z_0\in \C\!\!\!\!&:&  \textrm{the IMEX GLM is stable for fixed} \\
& & z_1=-|y|/\tan(\alpha)+iy
\end{array}
\right\}.
\end{equation}
Observe that 
\begin{equation} \label{eq2.15}
\mathcal{S}_{\alpha}=\bigcap_{y\in\R}\mathcal{S}_{\alpha,y}.
\end{equation}
Observe also that the region $\mathcal{S}_{\alpha,0}$ is independent of
$\alpha$, and corresponds to the region of absolute stability
of the explicit method with coefficients $\mbf{c}$, $\mbf{A}$, $\mbf{U}$, $\mbf{B}$, and $\mbf{V}$.
This region will be denoted by $\mathcal{S}_E$.
We have
\begin{equation} \label{eq2.16}
\mathcal{S}_{\alpha}\subset \mathcal{S}_E,
\end{equation} 
and we will look for IMEX DIMSIMs for which the stability region $\mathcal{S}_{\alpha}$ contains a 
large part of the stability region $\mathcal{S}_E$ of the explicit method.

All these regions $\mathcal{S}_E$, $\mathcal{S}_{\alpha,y}$, and
$\mathcal{S}_{\alpha}$, for fixed $y\in\R$ and $\alpha\in(0,\pi/2]$, can
be determined by the algorithms developed in a recent paper \cite{cjsz14a}.
These algorithms are based on some variants of boundary locus method
to compute the boundaries $\partial\mathcal{S}_E$, $\partial\mathcal{S}_{\alpha,y}$,
and $\partial\mathcal{S}_{\alpha}$, of the regions
$\mathcal{S}_E$, $\mathcal{S}_{\alpha,y}$,
and $\mathcal{S}_{\alpha}$.
We refer to the paper \cite{cjsz14a} for a detailed description of these 
algorithms.
The areas of $\mathcal{S}_E$ and $\mathcal{S}_{\alpha}$ can be computed
by numerical integration in polar coordinates. We refer again to
\cite{cjsz14a} for a detailed description of this process.

\setcounter{equation}{0}
\setcounter{figure}{0}
\setcounter{table}{0}

\section{Transformed IMEX DIMSIMs} \label{sec3}
Similarly as in \cite{cij18,ij17},
to increase our chances of finding SSP explicit GLMs with large
SSP coefficients $\mathcal{C}$ we consider a very general class
of transformed IMEX methods. These schemes are defined by
multiplying the relation for $y^{[n+1]}$ in (\ref{eq1.5}) by
$\mbf{T}\otimes \mbf{I}$, where $\mbf{T}\in\R^{r\times r}$, and $\det(\mbf{T})\neq 0$.
This leads to
\begin{equation} \label{eq3.1}
\begin{array}{l}
Y^{[n+1]}=h(\mbf{A}\otimes \mbf{I})f\big(Y^{[n+1]}\big)+h(\mbf{A}^{*}\otimes \mbf{I})g\big(Y^{[n+1]}\big) \\
\quad\quad\quad \ + \ (\mbf{U}\otimes \mbf{I})(\mbf{T}^{-1}\otimes \mbf{I})(\mbf{T}\otimes \mbf{I})y^{[n]}, \\ [3mm]
(\mbf{T}\otimes \mbf{I})y^{[n+1]}=h(\mbf{T}\otimes \mbf{I})(\mbf{B}\otimes \mbf{I})f\big(Y^{[n+1]}\big)
+h(\mbf{T}\otimes \mbf{I})(\mbf{B}^{*}\otimes \mbf{I})g\big(Y^{[n+1]}\big) \\
\quad\quad\quad \ + \ (\mbf{T}\otimes \mbf{I})(\mbf{V}\otimes \mbf{I})
(\mbf{T}^{-1}\otimes \mbf{I})(\mbf{T}\otimes \mbf{I})y^{[n]},
\end{array}
\end{equation}
$n=0,1,\ldots,N-1$.
Putting 
$$
\bar{y}^{[n+1]}=(\mbf{T}\otimes \mbf{I})y^{[n+1]},
\quad
\bar{y}^{[n]}=(\mbf{T}\otimes \mbf{I})y^{[n]},
$$
the equation (\ref{eq3.1}0 can be written in the form
\begin{equation} \label{eq3.2}
\begin{array}{l}
Y^{[n+1]}=h(\bar{\mbf{A}}\otimes \mbf{I})f(Y^{[n+1]})+h(\bar{\mbf{A}^{*}}\otimes \mbf{I})g(Y^{[n+1]})
+(\bar{\mbf{U}}\otimes \mbf{I})\bar{y}^{[n]}, \\
\bar{y}^{[n+1]}=h(\bar{\mbf{B}}\otimes \mbf{I})f(Y^{[n+1]})+h(\bar{\mbf{B}^{*}}\otimes \mbf{I})g(Y^{[n+1]})
+(\bar{\mbf{V}}\otimes \mbf{I})\bar{y}^{[n]},
\end{array}
\end{equation}
where the transformed coefficient matrices 
$\bar{\mbf{A}}$, $\bar{\mbf{A}^{*}}$, $\bar{\mbf{U}}$, 
$\bar{\mbf{B}}$, $\bar{\mbf{B}^{*}}$, and $\bar{\mbf{V}}$, 
are defined by
\begin{equation}\label{eq3.3}
\begin{split}
\bar{\mbf{A}}=\mbf{A},
\quad
\bar{\mbf{A}^{*}} & =\mbf{A}^{*},
\quad
\bar{\mbf{U}}=\mbf{U}\mbf{T}^{-1},
\\
\bar{\mbf{B}}=\mbf{T}\mbf{B},
\quad
\bar{\mbf{B}^{*}} & =\mbf{T}\mbf{B}^{*},
\quad
\bar{\mbf{V}}=\mbf{T}\mbf{V}\mbf{T}^{-1}.
\end{split}
\end{equation}
It was demonstrated in \cite{cij18} that transformed explicit 
and implicit methods preserve the order $p$ and stage order
$q=p$ of the original schemes. As a result, it follows from \cite{zsb14},
that the transformed IMEX method (\ref{eq3.2}) preserve
the order $p$ and stage order $q=p$ of the original 
IMEX method (\ref{eq1.5}).

Transformed IMEX GLMs (\ref{eq3.2}) preserve also absolute stability properties
of the original IMEX schemes (\ref{eq1.5}). This follows from
$$
\begin{array}{lcl}
\bar{\mbf{M}}(z_0,z_1) & = & \bar{\mbf{V}}+(z_0\bar{\mbf{B}}+z_1\bar{\mbf{B}^{*}})
(\mbf{I}-z_0\bar{\mbf{A}}+z_1\bar{\mbf{A}^{*}})^{-1}\bar{\mbf{U}} \\
& = & 
\mbf{T}\mbf{V}\mbf{T}^{-1}+(z_0\mbf{T}\mbf{B}+z_1\mbf{T}\mbf{B}^{*})
(\mbf{I}-z_0\mbf{A}+z_1\mbf{A}^{*})^{-1}\mbf{U}\mbf{T}^{-1} \\
& = &
\mbf{T}\big(\mbf{V}+(z_0\mbf{B}+z_1\mbf{B}^{*})
(\mbf{I}-z_0\mbf{A}+z_1\mbf{A}^{*})^{-1}\mbf{U}\big)\mbf{T}^{-1} \\
& = &
\mbf{M}(z_0,z_1),
\end{array}
$$
which shows that the stability matrix $\bar{\mbf{M}}(z_0,z_1)$ of the transformed
method is similar to the stability matrix $\mbf{M}(z_0,z_1)$ of the 
original method.
Hence, it follows that
$$
\bar{p}(z_0,z_1)=\det\big(w\mbf{I}-\bar{\mbf{M}}(z_0,z_1)\big)
=\det\big(w\mbf{I}-\mbf{M}(z_0,z_1)\big)
=p(z_0,z_1),
$$
and we can conclude that the transformed explicit, implicit, and IMEX methods
have identical absolute stability properties as the original explicit, implicit, and 
the IMEX methods. 
However, SSP properties of the transformed explicit GLMs are, in general,
different from SSP properties of the original explicit methods, and we will
search for transformed explicit DIMSIMs with maximal 
SSP coefficients.
In addition, 
we will monitor the size of the region of absolute stability 
$\mathcal{S}_{\alpha}$ for $\alpha\in(0,\pi/2)$, preferably for $\alpha=\pi/2$,
of the IMEX schemes, assuming that the 
implicit part of the method is $A$-, or $L$-stable.

\setcounter{equation}{0}
\setcounter{figure}{0}
\setcounter{table}{0}

\section{Construction of SSP transformed IMEX DIMSIMs} \label{sec4}
In this section we investigate transformed SSP IMEX DIMSIMs 
of order $p=1,2,3$, and $4$, with $q=r=s=p$. 
Our aim is to construct IMEX methods whose explicit part has large SSP coefficient,
the implicit part is $A$- or $L$-stable, and the overall IMEX scheme
has large region of absolute stability.
These methods will be compared with transformed SSP DIMSIMs investigated
recently in \cite{ij17}.

For many examples of DIMSIMs constructed in the literature on the subject, the abscissa vector
$\mbf{c}$ has components uniformly distributed in the interval $[0,1]$, i.e.,
$$
\mbf{c}=\left[
\begin{array}{ccccc}
 0 & \frac{1}{s-1} & \cdots & \frac{s-2}{s-1} & 1 
\end{array}
\right]^T\in\R^s.
$$
In our search for IMEX schemes with good stability properties we relax this
condition and consider methods with abscissa vector $\mbf{c}$ of the more
general form with abscissas satisfying the condition
$$
0<c_1<c_2<\cdots<c_{s-1}<c_s=1.
$$
Then the last stage $Y^{[n]}_s$ of the method (\ref{eq1.5}) 
approximates the solution $y$ to (\ref{eq1.1}) at the point $t_n$. 
This simplifies the implementation of these methods since no special finishing procedure
is needed. This is discussed in more detail in \cite{ij17}.
However, these methods still need starting procedures to compute
sufficiently accurate starting vector $y^{[0]}$.
Starting procedures for GLMs are discussed in \cite{cij17,cij18,ij15,ij17} 
and for IMEX methods in \cite{bcjp18}.

The case $p=1$ is not very interesting and it is not properly allowed for this class of methods 
because it is not possible to construct an IMEX DIMSIMs 
with $p=q=r=s=1$ where the implicit and the explicit 
parts share the same abscissa vector.
However, if this last condition is relaxed (that is 
the implicit and the explicit part are allowed to have 
different abscissa vectors), 
%
then the explicit part is the forward Euler method, while the 
implicit part is $A$-stable for $\lambda\geq 1/2$ and $L$-stable 
for $\lambda=1$. This last choice corresponds to the backward Euler method.
In this case, the stability function $R(z_0,z_1)$ is the product of the 
stability function of the explicit method and the stability function 
of the implicit one. For this reason, the stability region $\mathcal{S}_{\frac{\pi}{2}}$
is equal to the stability region $\mathcal{S}_E$ of the explicit method, 
which corresponds to the forward Euler method for any $\lambda$.
The IMEX scheme corresponding to $\lambda=1/2$ will be denoted
by IMEX DIMSIM1A, and to $\lambda=1$ by IMEX DIMSIM1L. 

For order $p=2$ and $p=3$
we succeeded in obtaining methods with
SSP explicit part and SSP coefficients close to that obtained in \cite{ij17},
and with $A$- or $L$-stable implicit part.
Unfortunately, in some cases these IMEX methods have quite small stability region  
$\mathcal{S}_{\frac{\pi}{2}}$ with respect to the stability region 
$\mathcal{S}_E$ of the explicit part.
However, larger $\mathcal{S}_{\frac{\pi}{2}}$ stability regions can be obtained
if one is willing to accept smaller SSP coefficients.

We have searched for IMEX schemes with large SSP coefficients of the explicit
part and $A$- or $L$-stability of the implicit part by solving the
minimization problem (\ref{eq2.6}) subject to the nonlinear constrains
(\ref{eq2.4}), and the constrains (\ref{eq2.8}) required for
$A$-stability, or the requirement that
the polynomials $p_1(z),p_2(z),\ldots,p_s(z)$ appearing in (\ref{eq2.7}) 
have degree less than $s$,
which is required for $L$-stability.
These minimization problems were solved using the MATLAB function
\texttt{fmincon} with randomly generated 
initial guesses.

The results of our numerical searches are summarized in Tables~\ref{tb4.1} and \ref{tb4.2},
where we have listed SSP coefficients $\mathcal{C}$, efficient SSP coefficients
$\mathcal{C}_{eff}$, $\textrm{area}(\mathcal{S}_E)$, $\textrm{area}(\mathcal{S}_{\pi/2})$,
and intervals of absolute stability $\textrm{int}(\mathcal{S}_E)$, and
$\textrm{int}(\mathcal{S}_{\pi/2})$. These tables correspond to methods which 
have a good balance between  
area of the stability region $\mathcal{S}_{\frac{\pi}{2}}$ and magnitude
of SSP coefficient. Table~\ref{tb4.1} corresponds to IMEX schemes for which 
the implicit part is $A$-stable, and Table~\ref{tb4.2} corresponds to IMEX schemes 
for which the implicit part is $L$-stable. The corresponding methods of order $p=2$
$p=3$, and $p=4$
are denoted by IMEX DIMSIM2A, IMEX DIMSIM2L, IMEX DIMSIM3A, IMEX DIMSIM3L, 
and IMEX DIMSIM4A. 
The coefficients of these methods are listed in the Appendix.

\begin{table}
$$
\hspace{-1.7cm}
\begin{array}{|c|c|c|c|c|c|c|c|}
\hline
\textrm{Method}  &  \phantom{++}\mathcal{C}\phantom{++} & \phantom{+}\mathcal{C}_{eff}\phantom{+}  & 
\textrm{area}(\mathcal{S}_E) & \textrm{area}(\mathcal{S}_{\pi/2}) &
\phantom{|}\textrm{int}(\mathcal{S}_E)\phantom{|} & \textrm{int}(\mathcal{S}_{\pi/2}) & 
\textrm{area}(\mathcal{S}_{RK}) \\
\hline 
\textrm{IMEX DIMSIM1A} & 1    & 1    &  3.14 & 3.14 &    (-2,0) &    (-2,0) & 3.14 \\
\textrm{IMEX DIMSIM2A} & 1.38 & 0.69 &  7.14 & 4.66 & (-2.87,0) & (-2.87,0) & 5.87 \\
\textrm{IMEX DIMSIM3A} & 0.99 & 0.33 &  9.68 & 2.18 & (-3.57,0) & (-1.32,0) & 9.12 \\
\textrm{IMEX DIMSIM4A} & 0.51 & 0.13 &  9.68 & 0.15 & (-3.01,0) & (-0.30,0) & 12.70 \\
\hline
\end{array}
$$
\caption{SSP coefficient $\mathcal{C}$,  effective SSP coefficients $\mathcal{C}_{eff}$,
$\textrm{area}(\mathcal{S_E})$ $\textrm{area}(\mathcal{S}_{\pi/2})$, 
$\textrm{int}(\mathcal{S}_E)$, $\textrm{int}(\mathcal{S}_{\pi/2})$, 
and $\textrm{area}(\mathcal{S}_{RK})$,
for transformed IMEX SSP DIMSIMs
with $p=q=r=s=1$, $p=q=r=s=2$, $p=q=r=s=3$, and $p=q=r=s=4$,
with $A$-stable implicit part.} \label{tb4.1}
\end{table}

\begin{table}
$$
\hspace{-1.7cm}
\begin{array}{|c|c|c|c|c|c|c|c|}
\hline
\textrm{Method}  &  \phantom{++}\mathcal{C}\phantom{++} & \phantom{+}\mathcal{C}_{eff}\phantom{+}  & 
\textrm{area}(\mathcal{S}_E) & \textrm{area}(\mathcal{S}_{\pi/2}) &
\phantom{|}\textrm{int}(\mathcal{S}_E)\phantom{|} & \textrm{int}(\mathcal{S}_{\pi/2}) & 
\textrm{area}(\mathcal{S}_{RK}) \\
\hline 
\textrm{IMEX DIMSIM1L} & 1    &    1  & 3.14 & 3.14 &    (-2,0) &    (-2,0) & 3.14 \\
\textrm{IMEX DIMSIM2L} & 1.17 & 0.59  & 7.46 & 7.34 & (-3.01,0) & (-3.01,0) & 5.87 \\
\textrm{IMEX DIMSIM3L} & 0.85 & 0.28  & 9.52 & 3.84 & (-4.10,0) & (-1.85,0) & 9.12 \\
\hline
\end{array}
$$
\caption{SSP coefficient $\mathcal{C}$,  effective SSP coefficients $\mathcal{C}_{eff}$,
$\textrm{area}(\mathcal{S_E})$ $\textrm{area}(\mathcal{S}_{\pi/2})$, 
$\textrm{int}(\mathcal{S}_E)$, $\textrm{int}(\mathcal{S}_{\pi/2})$, 
and $\textrm{area}(\mathcal{S}_{RK})$,
for transformed IMEX SSP DIMSIMs
with $p=q=r=s=1$, $p=q=r=s=2$, $p=q=r=s=3$, and $p=q=r=s=4$,
with $L$-stable implicit part.} \label{tb4.2}
\end{table}

\begin{figure}[t!h!b!]
\begin{center}
\scalebox{0.52}{\includegraphics{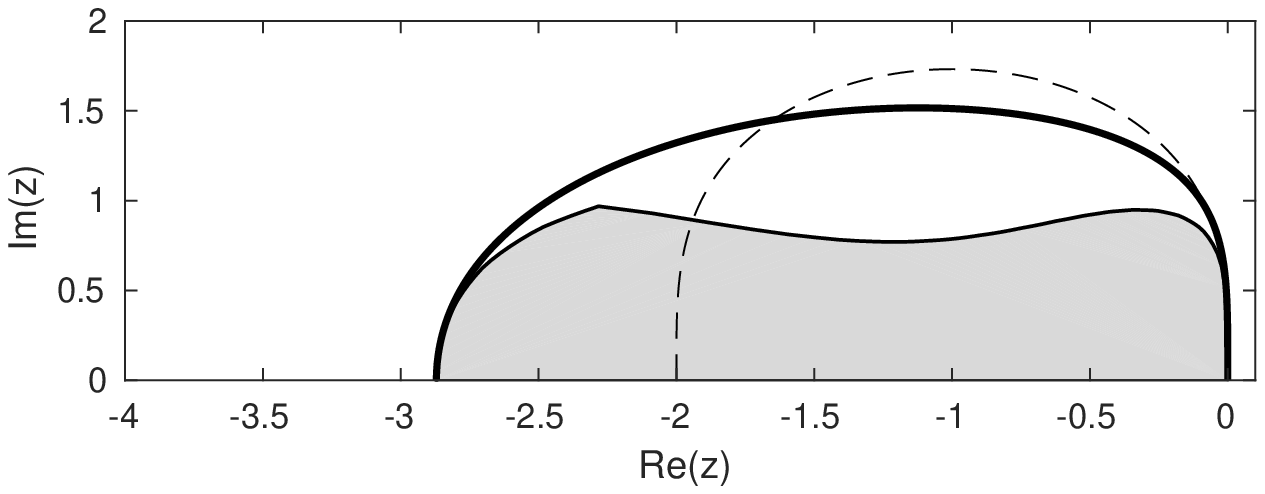}}
\scalebox{0.52}{\includegraphics{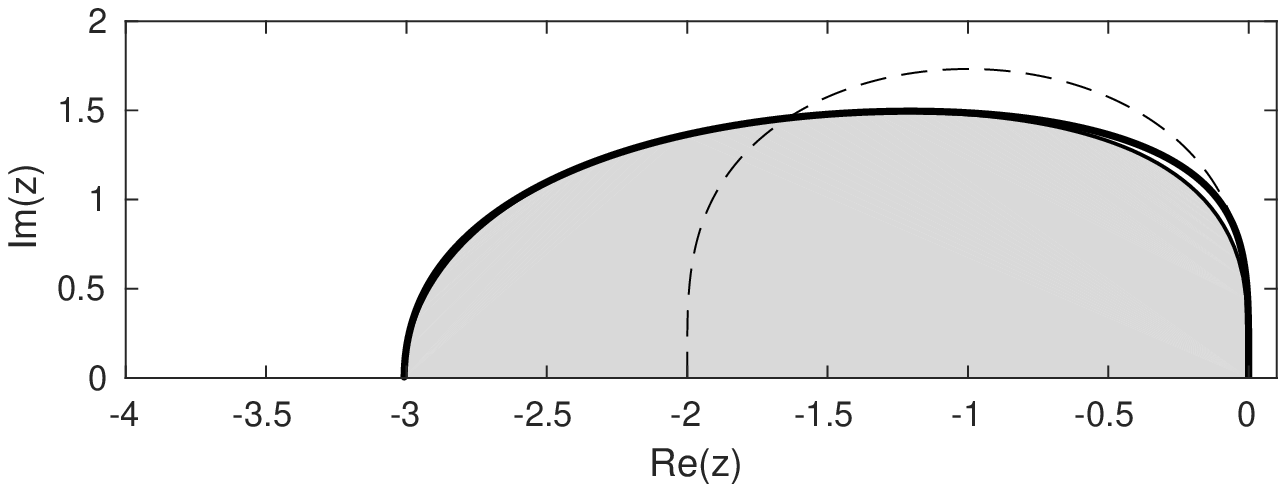}}
\caption{{\bf Left:} Stability region $\mathcal{S}_E$ (thick line), 
stability region $\mathcal{S}_{\pi/2}$ (shaded region),
of IMEX DIMSIM2A, and stability region of RK method of
order $p=2$ (dashed line). 
{\bf Right:} Stability region $\mathcal{S}_E$ (thick line), 
stability region $\mathcal{S}_{\pi/2}$ (shaded region),
of IMEX DIMSIM2L, and stability region of RK method of
order $p=2$ (dashed line). 
} \label{fig4.1ab} 
\end{center}
\end{figure}

\begin{figure}[t!h!b!]
\begin{center}
\scalebox{0.52}{\includegraphics{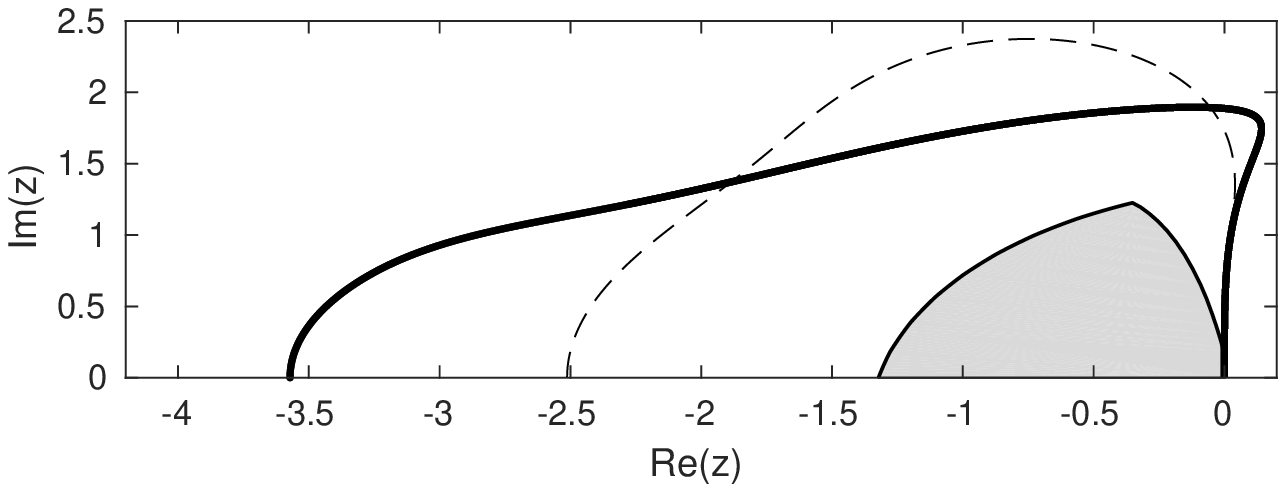}}
\scalebox{0.52}{\includegraphics{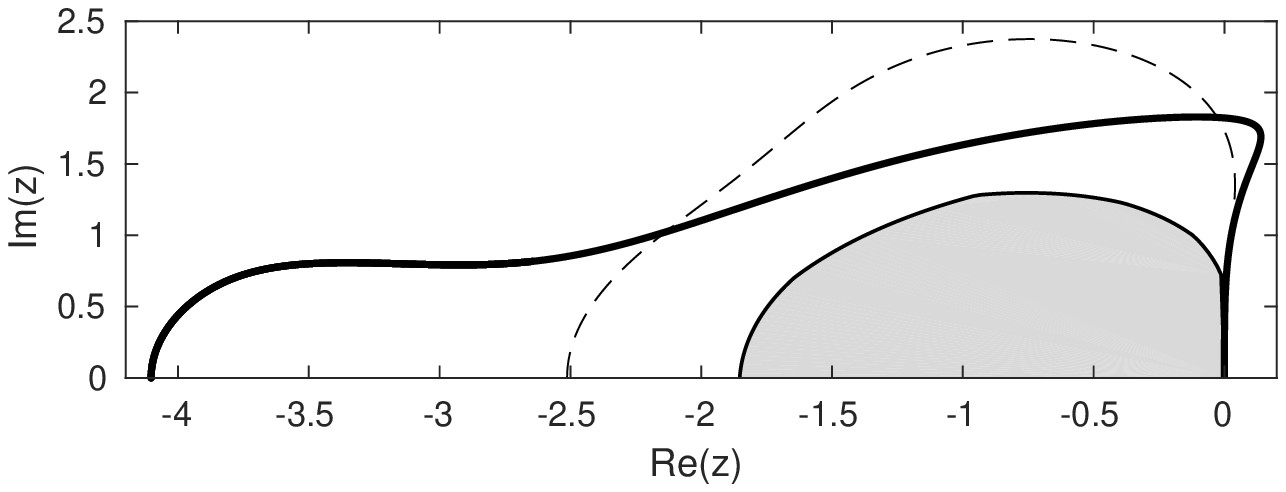}}
\caption{{\bf Left:} Stability region $\mathcal{S}_E$ (thick line), 
stability region $\mathcal{S}_{\pi/2}$ (shaded region),
of IMEX DIMSIM3A, and stability region of RK method of
order $p=3$ (dashed line). 
{\bf Right:} Stability region $\mathcal{S}_E$ (thick line), 
stability region $\mathcal{S}_{\pi/2}$ (shaded region),
of IMEX DIMSIM3L, and stability region of RK method of
order $p=3$ (dashed line). 
} \label{fig4.2ab} 
\end{center}
\end{figure}

\begin{figure}[t!h!b!]
\begin{center}
\scalebox{0.52}{\includegraphics{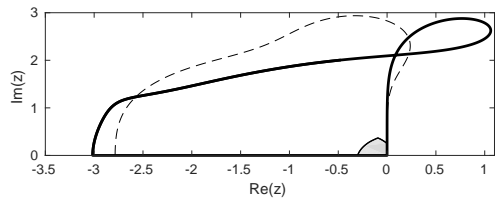}}
\caption{Stability region $\mathcal{S}_E$ (thick line), 
stability region $\mathcal{S}_{\pi/2}$ (shaded region),
of IMEX DIMSIM4A, and stability region of RK method of
order $p=4$ (dashed line). 
} \label{fig4.3} 
\end{center}
\end{figure}

The stability regions of these methods for order $p=2$, $3$, and $4$ are
reported in Figures~\ref{fig4.1ab}-\ref{fig4.3}.


\setcounter{equation}{0}
\setcounter{figure}{0}
\setcounter{table}{0}

\section{Numerical experiments} \label{sec5}
It has been shown in \cite{hr07} (but also \cite{ij17rk}) that IMEX RK can suffer 
from order reduction when applied to stiff problem. 
In order to confirm the good performances of the proposed methods
when applied to stiff problems, we solved several problems from literature,
such as shallow water equation \cite{ij17rk,pr05} with $\varepsilon=10^{-8}$, Schnakenberg reaction-diffusion 
\cite{bij17,hv03,sch79}, Van der Pol oscillator \cite{bij17,ij17rk} with $\varepsilon=10^{-6}$.  
In each considered case it has been confirmed that DIMSIM2A, DIMSIM2L, DIMSIM3A and DIMSIM3L
converge and achieve the expected order of convergence, 
while order reduction can occur for IMEX RK of order $p=2,3$ and $4$.
We also noticed that the performances of the method DIMSIM4A were good too, but
in some case, the behavior of the method was somehow erratic. This was probably 
motivated by the fact that the performances of this last method were sensitive 
to perturbation in the starting procedure, were IMEX RK methods sometimes were 
not enough to get good starting values.

For the sake of brevity, we report here detailed results for numerical resolution of
an advection-reaction problem, an adsorption-desorption problem, and  
a shallow water problem.

\subsection{Problem 1: advection-reaction} 
Consider next the linear advection-reaction equation \cite{bij17,cjsz14a,hv03}
\begin{equation} \label{eq9.2}
\left\{
\begin{array}{l} 
\ds\frac{\partial u}{\partial t}+\alpha_1\,\ds\frac{\partial u}{\partial x}
=-k_1u+k_2v+s_1, \\ [3mm]
\ds\frac{\partial v}{\partial t}+\alpha_2\,\ds\frac{\partial v}{\partial x}
=k_1u-k_2v+s_2, 
\end{array}
\right.
\end{equation}
$0\leq x\leq 1$, $0\leq t\leq 1$, with parameters
$$
\alpha_1=1,
\quad\alpha_2=0,
\quad
k_1=10^6,
\quad
k_2=2k_1,
\quad
s_1=0,
\quad
s_2=1,
$$
and with initial and boundary values
$$
u(x,0)=1+s_2x,
\quad
v(x,0)=\ds\frac{k_1}{k_2}u(x,0)+\ds\frac{s_2}{k_2},
\quad
0\leq x\leq 1,
$$
$$
u(0,t)=\gamma_1(t),
\quad
v(0,t)=\gamma_2(t),
\quad
0\leq t\leq 1.
$$
(Observe that the condition $v(0,t)=\gamma_2(t)$ does not have to be specified since
$\alpha_2=0$).
Discretization of (\ref{eq9.2}) in space variable $x$ on the uniform grid
\mbox{$x_i=i\Delta x$}, $i=0,1,\ldots,N$, $N\Delta x=1$,
leads to the initial value problem for the system
of ODEs of dimension $2N$, with non-stiff part
corresponding to the advection terms, and stiff part corresponding
to the reaction terms.

We consider the spatial discretization of (\ref{eq9.2}) which corresponds to the time dependent
Dirichlet data $\gamma_1(t)=1-\sin(12t)^4$ at the left boundary,
where $u_x$ is approximated by fourth-order central differences in the interior
domain and third-order finite differences at the boundary, as in \cite{bij17,hv03}.

\begin{figure}[t!h!b!]
\begin{center}
\scalebox{1}{\includegraphics{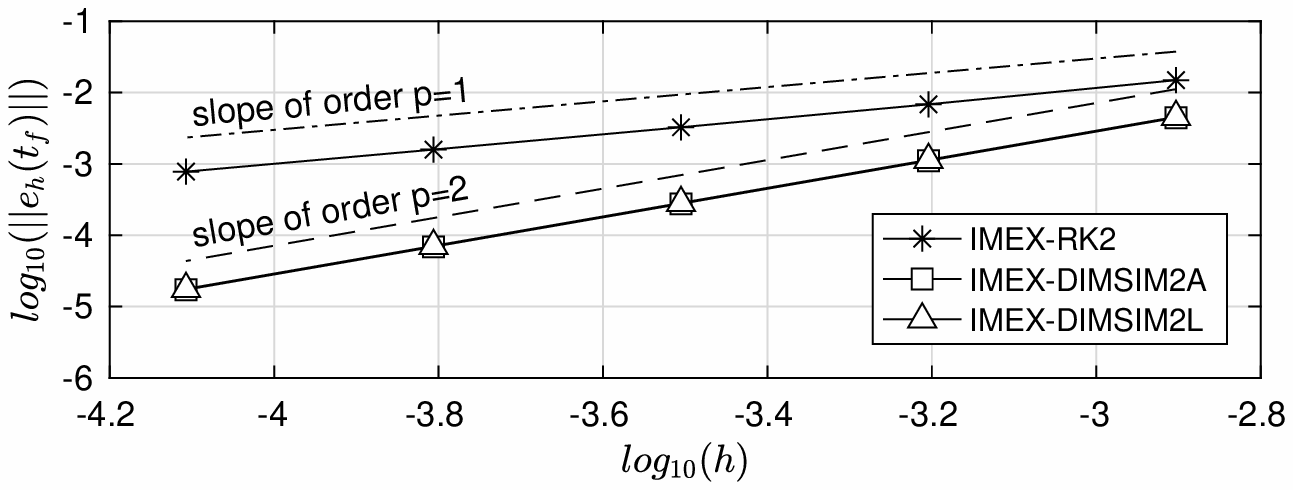}}\\ \vspace{2mm} 
\scalebox{0.98}{\includegraphics{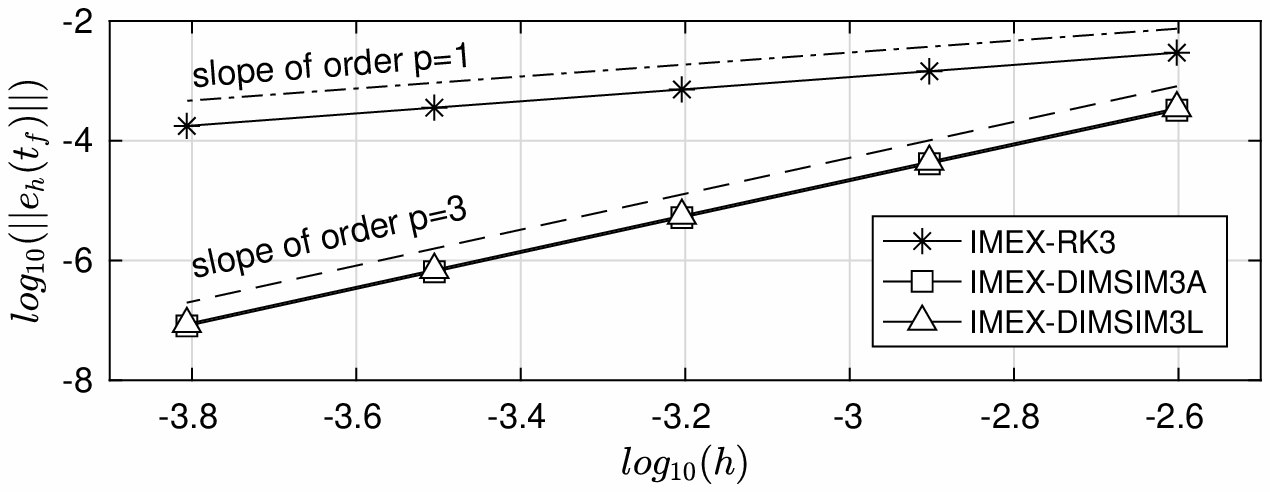}}\\ \vspace{2mm} 
\scalebox{1}{\includegraphics{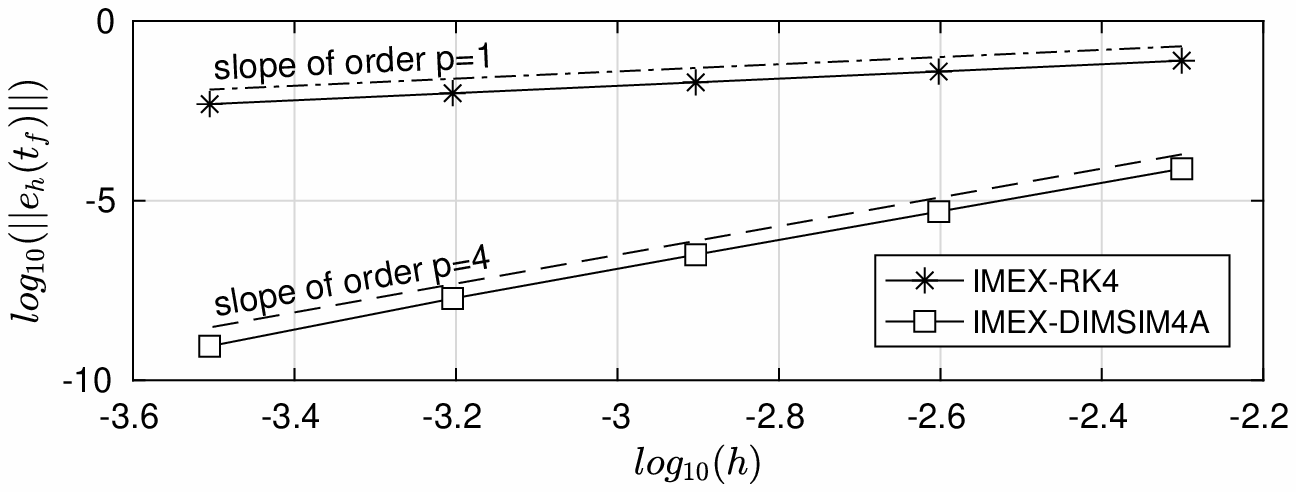}}
\caption{Error versus stepsize (double logarithmic scale plot) for 
SSP transformed IMEX DIMSIMs and IMEX RK,
applied to the discretization of the Advection-Reaction problem
(\ref{eq9.2}), with $N=401$, by fourth-order central differences
in the interior domain and by third-order finite differences at the boundaries. 
} \label{fig:5.1} 
\end{center}
\end{figure}

\begin{figure}[t!h!b!]
\begin{center}
\scalebox{1}{\includegraphics{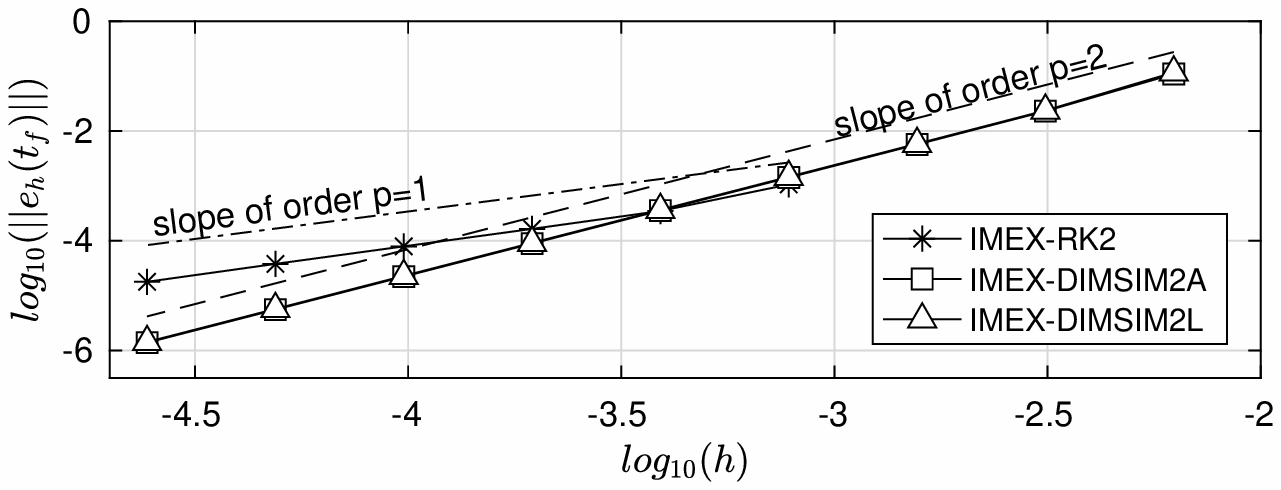}}\\ \vspace{2mm}
\scalebox{0.98}{\includegraphics{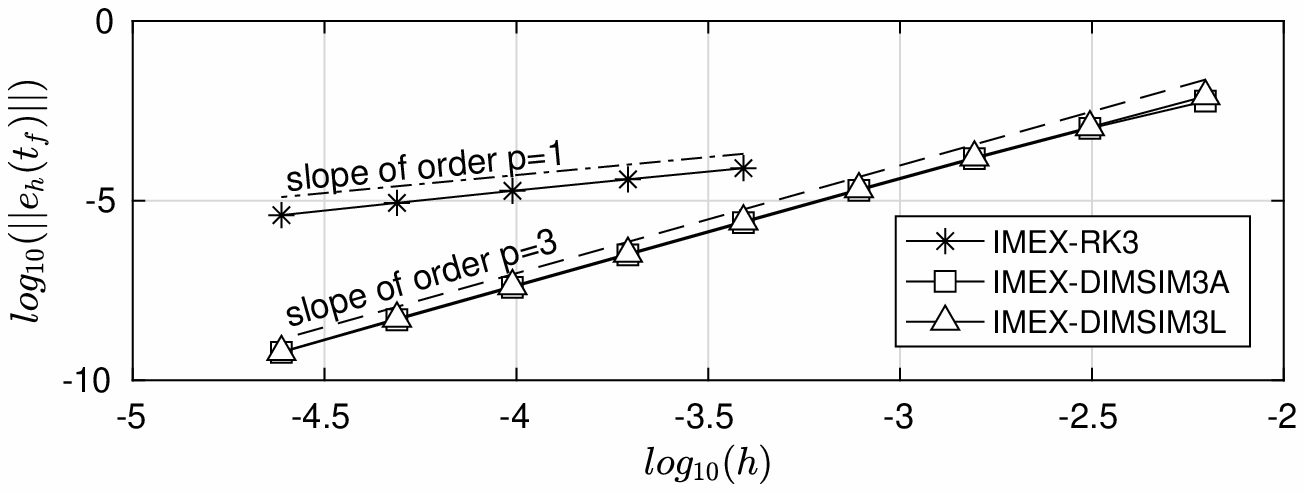}}\\ \vspace{2mm}
\scalebox{1}{\includegraphics{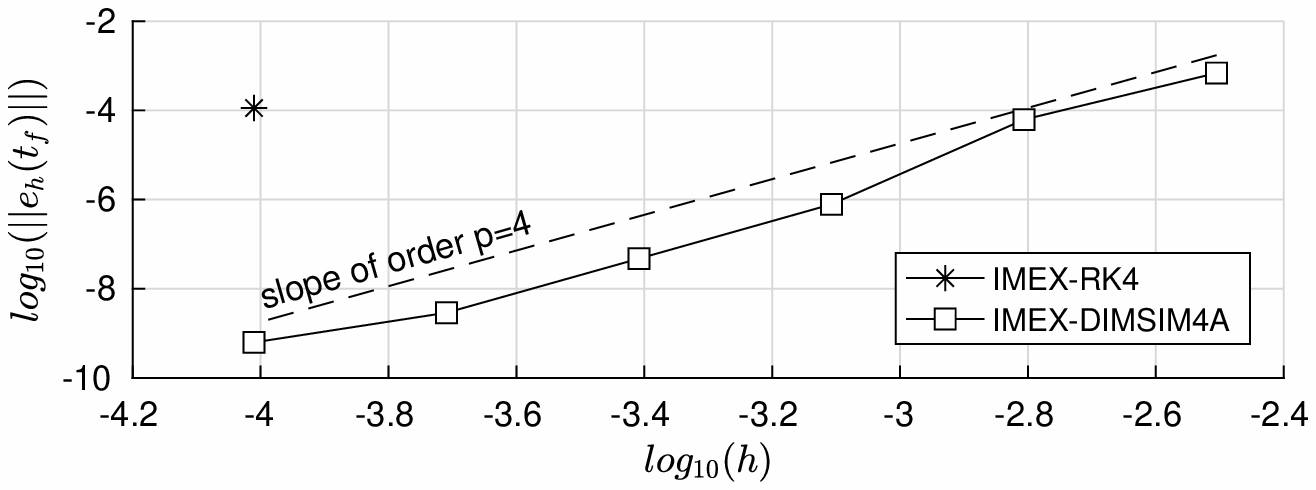}}
\caption{Error versus stepsize (double logarithmic scale plot) for 
SSP transformed IMEX DIMSIMs and IMEX RK,
applied to the discretization of the Adsorption-desorption problem
(\ref{eq9.5}), with $N=101$, by a WENO5 space discretization scheme. 
} \label{fig:5.2} 
\end{center}
\end{figure}

\begin{figure}[t!h!b!]
\begin{center}
\scalebox{1}{\includegraphics{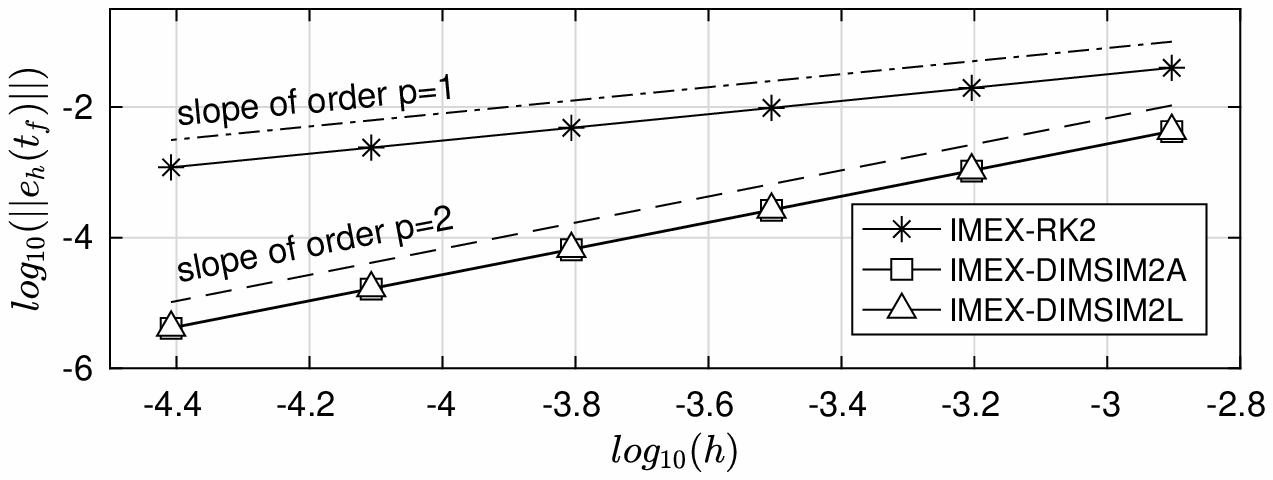}}\\ \vspace{2mm} 
\scalebox{0.98}{\includegraphics{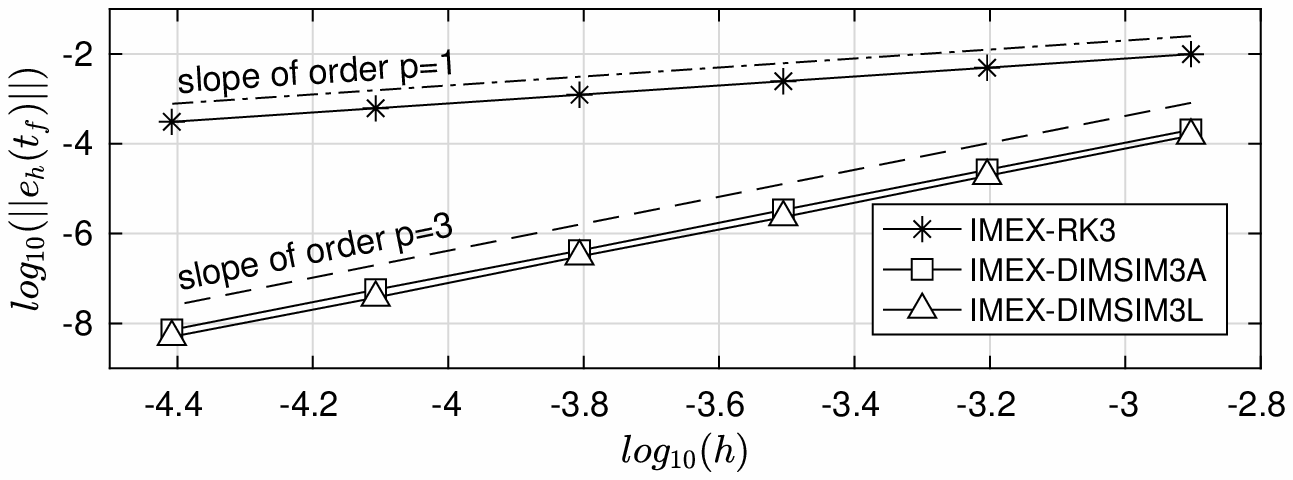}}\\ \vspace{2mm} 
\scalebox{1}{\includegraphics{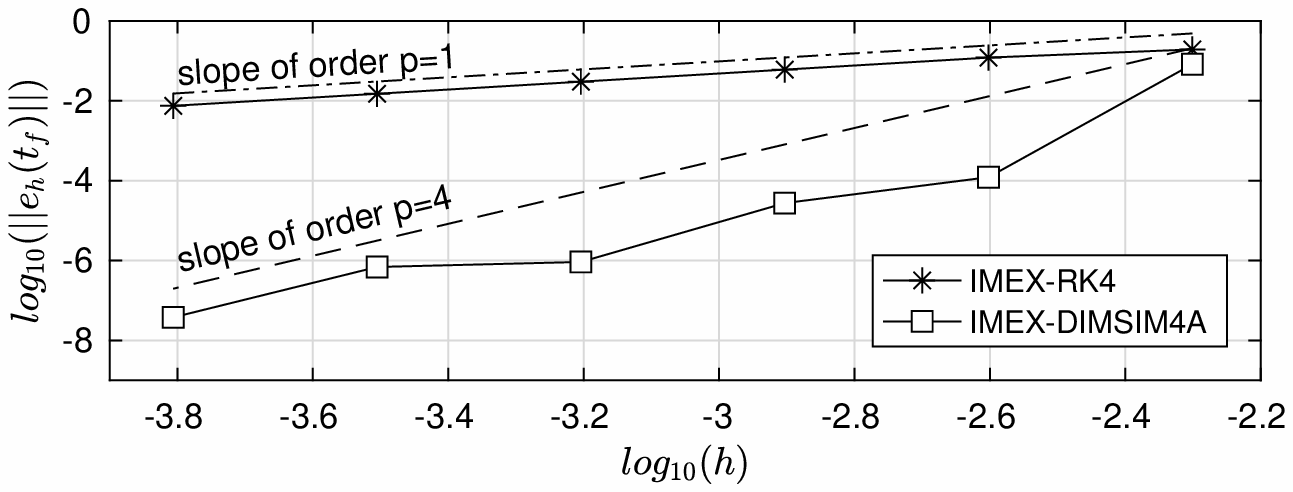}}
\caption{Error versus stepsize (double logarithmic scale plot) for 
SSP transformed IMEX DIMSIMs and IMEX RK,
applied to the discretization of the shallow water problem (\ref{eq:swe_1d_pr}), 
with $N=201$, by a WENO5 space discretization scheme, and $\epsilon=10^{-8}$.  
} \label{fig:6.1} 
\end{center}
\end{figure}

The numerical results for the discretization of (\ref{eq9.2}) with $N=400$ spatial 
points are presented in Figure~\ref{fig:5.1}, where for the sake of comparison
we also report the results obtained by IMEX RK methods constructed 
in \cite{ij17rk}.
For these tests, the reference solution was
computed by \texttt{ODEPACK} routine \texttt{DLSODAR} (\cite{h83}) 
with absolute tolerance and relative tolerance equal to $10^{-14}$ and
$10^{-12}$, respectively.

To start the integration, we used the starting procedure described in Section 2 of \cite{cij17},
where the required starting values have been computed, for methods of order
$p=2$ and $p=3$, by IMEX RK of the same order $p$, applied with a
suitable stepsize, for the method of order $p=4$, by \texttt{DLSODAR} 
with absolute tolerance and relative tolerance equal to $10^{-14}$ and
$10^{-12}$, respectively.

Figure~\ref{fig:5.1} shows that all the presented IMEX DIMSIMs  
achieve the expected order of convergence for this stiff system of ODEs
while order reduction to $p=1$ occurs for IMEX RK of the same order.

\subsection{Problem 2: adsorption-desorption model}
Following Hundsdorfer and Ruuth \cite{hr07} (see also \cite{hv03}) we consider
next the adsorption-desorption problem given by the equations
\begin{equation} \label{eq9.5}
\left\{
\begin{array}{l}
u_t+a(t)u_x=\kappa\big(v-\phi(u)\big), \\
v_t=-\kappa\big(v-\phi(u)\big),
\end{array}
\right.
\end{equation}
$0\leq x\leq 1$, $t\in[0,t_{end}]$, $t_{end}=1.25$, where
$\phi(u)=k_1u/(1+k_2u)$.
The initial values are
$
u(x,0)=v(x,0)=0,
$
$
0\leq x\leq 1,
$
and the boundary values are 
$$
\left\{ 
\begin{array}{lll}
u(0,t)=1-\cos^2(6\pi t), \quad & a\geq 0, \\
u(1,t)=0, & a<0.
\end{array}
\right.
$$
As in \cite{hr07} we choose the parameters $\kappa=10^6$,
$k_1=50$, $k_2=100$, and the velocity
$
a=a(t)=-\arctan\big(100(t-1)\big)/\pi.
$
Then $a(t)>0$ for $0\leq t\leq 1$, which corresponds
to the adsorption phase, and $a(t)<0$ for $t>1$, which
corresponds to the desorption phase.

As in \cite{bij17}, for the spatial discretization of $u_x$ we have
implemented the WENO5 scheme \cite{shu99} following the presentation
in \cite{ws07}. Further details can be found in \cite{bij17}.

The results of these tests are reported in 
Figure \ref{fig:5.2}  
where it is confirmed that all the presented IMEX DIMSIMs  
achieve the expected order of convergence for this stiff system of ODEs
while order reduction to $p=1$ occurs for IMEX RK of the same order. 
We also point out that several points on the line corresponding to IMEX RK 
of order $p=4$ are missing because it did not converge for several values 
of the stepsize $h$.

\subsection{Problem 3: shallow water model}
 We now consider a one-dimensional  model of shallow water flow (compare \cite{pr05,jin95}):
\begin{equation}\label{eq:swe_1d_pr}
\left\{
 \begin{array}{l}
  \ds \dt h +\dx(hv) = 0,  \vspace{3mm}  \\ 
  \ds \dt(hv) +\dx \left( h + \frac{1}{2}h^2 \right) = \ds \frac{1}{\varepsilon}\left( \frac{h^2}{2} - hv \right), \\
 \end{array}
\right.
\end{equation} 
where $h$ is the water height with respect to the bottom and $hv$ is the flux 
of the velocity field.
We use periodic boundary conditions and initial conditions
at $t_0 = 0$
\begin{equation}\label{eq:swe_bc}
h(0 , x) = 1 + \frac{1}{5}sin(8\pi x), \quad hv(0 , x) = \frac{1}{2} h(0 , x)^2 ,
\end{equation}
with $x\in[0,1]$. For this problem the space derivative was discretized by  
a fifth order finite difference weighted essentially non-oscillatory (WENO) scheme following 
the implementation described in \cite{shu99}.

The numerical results for the discretization of (\ref{eq:swe_1d_pr}) with $\epsilon=10^{-8}$ 
and $N=201$ spatial 
points are presented in Figure~\ref{fig:6.1}, where for the sake of comparison
we also report the results obtained by IMEX Runge-Kutta methods constructed 
in \cite{ij17rk}.
For these tests, the reference solution was
computed by \texttt{ODEPACK} routine \texttt{DLSODAR} (\cite{h83}) 
with absolute tolerance and relative tolerance equal to $10^{-14}$ and
$10^{-12}$, respectively.

To start the integration, we used the starting procedure described in Section 2 of \cite{cij17},
where the required starting values have been computed, for methods of order
$p=2$ and $p=3$, by IMEX RK of the same order $p$, applied with a
suitable stepsize, for the method of order $p=4$, by \texttt{DLSODAR} 
with absolute tolerance and relative tolerance equal to $10^{-14}$ and
$10^{-12}$, respectively.

The results reported in Figure~\ref{fig:6.1} confirm that the IMEX RK
just has the asymptotic preserving property, while the methods proposed in this 
paper are also asymptotically accurate in the stiff limit 
for $\varepsilon \to 0$ (compare \cite{pr05}). In other words, 
the IMEX RK methods converge, but the order drop to $p=1$, while
all the presented IMEX DIMSIMs achieve the expected order of convergence 
and no order reduction occurs.

\appendix
\renewcommand{\thesection}{\Alph{section}.\arabic{section}}
\setcounter{section}{0}

\begin{appendices}

\section{Appendix}
In this Appendix we report the coefficients of the methods of order $p=2,3,4$, 
described in Section \ref{sec4}.

\subsection{Coefficients of method DIMSIM2A}
\begin{equation}\nonumber 
\begin{split}
 \mbf{c}= & \left[ 
\begin{array}{cc}
 0.5207015987954746 & 1 \\
\end{array}
\right]^T,
\\
 \bar{\mbf{A}}= & \left[
\begin{array}{cc}
 0 & 0 \\
 0.6335780271090006 & 0 \\
\end{array}
\right],
\\
 \bar{\mbf{A^*}}= & \left[
\begin{array}{cc}
 0.9756662942012514 & 0 \\
 1.065344873186484 & 0.9756662942012514 \\
\end{array}
\right],
\\
\bar{\mbf{U}}= & \left[
\begin{array}{cc}
 1 & 0 \\
 0.8760323181723925 & 1 \\
\end{array}
\right]^T,
\\
\bar{\mbf{V}}= & \left[
\begin{array}{cc}
 0.8035259425918053 & 1.584881273180670 \\
 0.09961124839144930 & 0.1964740574081947 \\
\end{array}
\right].
\end{split}
\end{equation}

\subsection{Coefficients of method DIMSIM2L}
\begin{equation}\nonumber 
\begin{split}
\mbf{c}= & \left[
\begin{array}{cc}
 0.5725000000000000 & 1 \\
\end{array}
\right]^T,\\
\bar{\mbf{A}}= & \left[
\begin{array}{cc}
 0 & 0 \\
 0.5507246376811594 & 0 \\
\end{array}
\right],\\
\bar{\mbf{A^*}}= & \left[
\begin{array}{cc}
 0.4025509997331064 & 0 \\
 0.3054637337141530 & 0.4025509997331064 \\
\end{array}
\right],\\
\bar{\mbf{U}}= & \left[
\begin{array}{cc}
 1 & 0 \\
 0.8970000000000000 & 1 \\
\end{array}
\right]^T,\\
\bar{\mbf{V}}= & \left[
\begin{array}{cc}
 0.7976747326679189 & 1.964322983806612 \\
 0.08216049746479565 & 0.2023252673320811 \\
\end{array}
\right].
\end{split}
\end{equation}

\subsection{Coefficients of method DIMSIM3A}
\begin{equation} \nonumber 
\begin{split}
  \mbf{c}= & \left[
  \begin{array}{ccc}
  0.3785922442536512 & 0.7369632894601272 & 1 \\
  \end{array}
  \right]^T,\\
  \bar{\mbf{A}}= & \left[
  \begin{array}{ccc}
  0 & 0 & 0 \\
  0.6105030326964779 & 0 & 0 \\
  0.5054775907409634 & 0.3826213150653439 & 0 \\
  \end{array}
  \right],\\
  \bar{\mbf{A^*}}= & \left[
  \begin{array}{ccc}
  0.5023463944444552 & 0 & 0 \\
  -0.8899211224523407 & 0.5023463944444552 & 0 \\
  -3.305290943287502 & 0.4193402392399124 & 0.5023463944444552 \\
  \end{array}
  \right],\\
  \bar{\mbf{U}}= & \left[
  \begin{array}{ccc}
  1 & 0 & 0 \\
  0.6070215241878391 & 1 & 0 \\
  0.5361152778084712 & 1.091180739129647 & 1 \\
  \end{array}
  \right]^T,\\
  \bar{\mbf{V}}= & \left[
  \begin{array}{ccc}
  0.5418838673478645 & 0.9017144383487438 & 2.958352027358458 \\
  0.2129486962575630 & 0.3543543656001081 & 1.162568670627143 \\
  0.01900613148571312 & 0.03162689316015439 & 0.1037617670520274 \\
  \end{array}
  \right].
\end{split}
\end{equation}

\subsection{Coefficients of method DIMSIM3L}
\begin{equation} \nonumber 
 \begin{split}
\mbf{c}= & \left[
\begin{array}{ccc}
 0.4020684033460171 & 0.7554528159803609 & 1 \\
\end{array}
\right]^T,\\
\bar{\mbf{A}}= & \left[
\begin{array}{ccc}
 0 & 0 & 0 \\
 0.5925366351567699 & 0 & 0 \\
 0.5582112117594124 & 0.3256969821842126 & 0 \\
\end{array}
\right],\\
\bar{\mbf{A^*}}= & \left[
\begin{array}{ccc}
 0.5201730949739405 & 0 & 0 \\
 -1.082981144838764 & 0.5201730949739405 & 0 \\
 -2.860648399647160 & 0.2917933416909193 & 0.5201730949739405 \\
\end{array}
\right],\\
\bar{\mbf{U}}= & \left[
\begin{array}{ccc}
 1 & 0 & 0 \\
 0.6343850217261301 & 1 & 0 \\
 0.5123644514467803 & 1.138668063964801 & 1 \\
\end{array}
\right]^T,\\
\bar{\mbf{V}}= & \left[
\begin{array}{ccc}
 0.4816666646770200 & 0.7031253548332313 & 3.663136087971684 \\
 0.1761045471411361 & 0.2570731613311589 & 1.339297421217996 \\
 0.03435316450098294 & 0.05014791919551827 & 0.2612601739918211 \\
\end{array}
\right].
 \end{split}
\end{equation}

\subsection{Coefficients of method DIMSIM4A}
\begin{changemargin}{-0.5cm}{2cm}
\begin{equation} \nonumber 
\hspace*{-1.8cm}
\begin{split}
\mbf{c}= & \left[
\begin{array}{cccc}
 0.2561983471074380 & 0.4485981308411215 & 0.7622950819672131 & 1 \\
\end{array}
\right]^T,\\
\bar{\mbf{A}}= & \left[
\begin{array}{cccc}
 0 & 0 & 0 & 0 \\
 0.3245033112582781 & 0 & 0 & 0 \\
 0.1102941176470588 & 0.6486486486486486 & 0 & 0 \\
 0.3111111111111111 & 0.1603053435114504 & 0.4729729729729730 & 0 \\
\end{array}
\right],\\
\bar{\mbf{A^*}}= & \left[
\begin{array}{cccc}
 1.228571428571429 & 0 & 0 & 0 \\
 -2.659574468085106 & 1.228571428571429 & 0 & 0 \\
 -6.431818181818182 & -0.4444444444444444 & 1.228571428571429 & 0 \\
 -5.931034482758621 & -4.906250000000000 & 1.103448275862069 & 1.228571428571429 \\
\end{array}
\right],\\
\bar{\mbf{U}}= & \left[
\begin{array}{cccc}
 1 & 0 & 0 & 0 \\
 0.7011494252873563 & 1 & 0 & 0 \\
 0.2363213391750847 & 0.3563218390804598 & 1 & 0 \\
 0.3704826947154125 & 0.5083355703606088 & 0.6222222222222222 & 1 \\
\end{array}
\right]^T,\\
\bar{\mbf{V}}= & \left[
\begin{array}{rrrr}
 0.3181770223788457 & 1.319227410800732 & 0.2619374293792898 & 1.680623378297797 \\
 0.09508738599827574 & 0.3942518698944718 & 0.07828015130875329 & 0.5022552624798014 \\
 0.2091032901032768 & 0.8669852710621154 & 0.1721430978104653 & 1.104491692074865 \\
 0.02185292729383308 & 0.09060673356209266 & 0.01799029847272758 & 0.1154280099162172 \\
\end{array}
\right].
 \end{split}
\end{equation}
\end{changemargin}

\end{appendices}

\end{document}